\documentstyle{amsppt}
\define\={\overset\text{def}\to=}
\define\C{{\Bbb C}}
\define\R{{\Bbb R}}

\define\re{{\operatorname{Re}}}
\define\im{{\operatorname{Im}}}
\define\unim{{\sqrt{-1}}}

\redefine\D{{\Cal D}}
\redefine\L{{\Cal L}}
\define\A{{\Cal A}}

\define\X{{\Cal X}}
\define\Y{{\Cal Y}}
\define\cZ{{\Cal Z}}
\redefine\H{{\Cal H}}
\define\derI#1#2#3{\nabla_{#1}#2 |_{#3}}
\define\derII#1#2#3{D_{#1} #2 |_{#3}}
\define\F#1{{\Cal F}^{(#1)}}
\define\doubledot#1{\dot v_{#1}}

\define\gl{\goth{gl_n(\C)}}
\redefine\u{\goth u_{n-1}\oplus \R}

\redefine\P{{\Bbb P}}
\define\TM{T^{o}M}
\define\PTM{{\Bbb P}TM}
\define\TTM{\widetilde{T}(\PTM)}
\define\CTTM{\A^0(\TTM)}

\define\UF{U_{F}(M)}
\define\UUF#1{U_{F}^{(#1)}(M)}

\define\bh{\text{h}}
\define\bH{\text{H}}
\define\gh{\text{\bf H}}

\TagsOnRight

\topmatter

\title
The   Hermitian connection
 and the  Jacobi fields\\
of a complex Finsler manifold
\endtitle

\author
A. Spiro
\endauthor

\address
\phantom{ }\newline
A. Spiro\newline
Dipartimento di Matematica e Fisica\newline
Universit\`a di Camerino\newline
Via Madonna delle Carceri\newline
62032 Camerino (Macerata)\newline
ITALY\newline
\phantom{ }\newline
\phantom{ }
\endaddress

\email
spiro\@campus.unicam.it
\endemail


\leftheadtext{A. Spiro}
\rightheadtext{The  Hermitian connection
of a complex Finsler manifold}

\abstract 
It is proved that all  invariant functions 
of a complex Finsler manifold can be totally recovered from the torsion and 
curvature of the connection introduced by Kobayashi
for holomorphic vector bundles with complex Finsler structures. 
Equations of the geodesics and
Jacobi fields of a generic complex Finsler manifold, 
expressed by means Kobayashi's connection, are also 
derived.
\endabstract

\endtopmatter

\document
\subhead 1. Introduction
\endsubhead
\bigskip
A complex Finsler manifold $(M, J, F)$ is a complex manifold 
$(M, J)$ endowed with a complex Finsler metric $F$, which is 
a continuous function 
$F\: TM \longrightarrow \R^+$ that is smooth on $TM\setminus\{\text{zero section}\}$ and 
verifies the following two properties:
\roster
\item"a)" 
 $F(u) >0$ for any $u>0$;
\item"b)" $F(\lambda u) = |\lambda|F(u)$
for any  $u \neq 0 $ and  $\lambda \in \C^*$.
\endroster
 In  this paper,  we will always
assume that $F$ is {\it strictly pseudoconvex\/}, i.e. 
that  any Finsler pseudo-sphere at a point $x$
$$S_x = \{ \ v\in T_xM\ :\ F(v) = 1\}$$
 is  strictly pseudoconvex as real hypersurface of $T_xM \simeq \C^n$.\par
The simplest examples of such Finsler manifold are  the Hermitian 
manifolds. In fact, 
if $g$ is an Hermitian metric on a complex manifold $(M, J)$,  the norm function
$$F_g: TM \to \R^+\ ,\qquad F(v) = \sqrt{g(v,v)}\tag 1.1$$
is a strictly pseudoconvex Finsler metric. In what follows, whenever 
 $F$ is as in (1.1), we will say that 
{\it $F$ is associated with the Hermitian metric $g$\/}.\par
Other important examples of complex Finsler manifolds are  
the  bounded convex domains in $\C^n$ with smooth boundary, endowed with their infinitesimal 
Kobayashi metric 
(see \cite{Le},  
\cite{Le1}, \cite{Pa}, \cite{Fa}, \cite{AP}, \cite{Sp}). Indeed,  the Kobayashi infinitesimal 
metric of any hyperbolic
complex manifold is a  "non-smooth" complex Finsler metric.
\medskip
In (\cite{Sp}), we introduced the concept of 
{\it adapted linear frames of a complex Finsler manifold  $(M, J, F)$\/}. 
We studied the properties of 
 the bundle $\pi: \UF\to M$ of all adapted frames  and we constructed
an absolute parallelism on $\UF$,  whose structure 
functions constitute a complete set of generators for 
the (local) invariants of $(M, J, F)$. \par
Such absolute parallelism consists of 
 a finite set of global vector fields 
$\{X_1, \dots,$ $X_{2n}$, $Y_1$, $\dots ,Y_p\}$
on $\UF$, which are preserved by any (local)
 diffeomorphism which is  lift of  a (local) biholomorphic isometry of $(M,J,F)$. 
It contains a subset $\{Y_1, \dots, Y_p\}$ of vector fields, which span the vertical distribution,
and another complementary subset 
$\{X_1, \dots , X_{2n}\}$, whose vector fields  
span the distribution $\H \subset T\UF$  of real subspaces
underlying the 
holomorphic distribution $T^{10}\UF \subset T^\C\UF$. By holomorphic distribution
$T^{10}\UF$  we use the standard meaning of 
 the set of subspaces $T^{10}_u\UF = T^\C_u\UF \cap T^{10}_uL^\C(M)$, 
where $T^{10}_uL^\C(M)$ is the subspace at $u\in \UF$ generated  by the holomorphic vector fields of 
the bundle of all complex linear frames $L^\C(M)$. \par
The distribution $\H$ is complementary to the 
vertical distribution. This allows
to interpret  any curve $u_t$ in $\UF$, which 
is tangent to  $\H_{\gamma_t}$
at any $t$,  as a 1-parameter family of adapted
frames, which represent a 'parallel transport'
 along the curve $\gamma = \pi\circ u: [a,b] \subset \R \to TM$. 
Using this  'parallel transport', 
one can define a covariant derivation
of  functions with values in $TM$ in the directions
of vectors tangent to $TM$. 
Note that, in case $F$ is associated 
with an Hermitian metric $g$, such covariant derivation 
reduces to the 
usual Hermitian covariant derivation of vector fields on $M$
in the direction of vectors tangent to $M$.\par
\par
\medskip
At the best of our knowledge, in the literature 
there exist other three different definitions
of Finslerian covariant derivation: namely, the one determined
by the absolute parallelism of
J.J. Faran in \cite{Fa}, the covariant derivation by M. Abate and G. Patrizio given in \cite{AP} 
and the one defined by S. Kobayashi in \cite{Ko}. 
Those definitions have the advantage to 
be defined using the expression of the Finsler metric in 
 complex coordinates  and  are therefore  suitable for explicit computations. 
On the other hand, our setting is established in a 
'totally coordinate-free language'\ and all objects we deal with (such as 
torsion and curvature)
immediately reduce to the corresponding  objects of Hermitian geometry, 
whenever  the Finsler metric is associated with an Hermitian metric.  \par
In this paper, we  show that 
our definition of Finslerian covariant derivation is 
strictly related with Kobayashi's definition and we derive the 
formulae which express the torsion and curvature of Kobayashi's connection 
in terms of the torsion and the
curvature 2-forms of the  absolute parallelism $\{ X_1, \dots, Y_p\}$
of $\UF$ (see \S 3 and \S 4). An important by-product of these 
formulae is a practical way to evaluate the structure functions of 
a complex Finsler manifold using coordinates: all computations reduce
to use the Kobayashi's expressions for the Finsler connection and 
Finsler curvature
given in \cite{Ko2} (see also \cite{Ko3}).\par
\medskip
We conclude  with  sections \S 5 and \S 6, in which 
 we  write  the equations of geodesics and of the Jacobi vector 
fields in terms of the Kobayashi's Finslerian connection and of its torsion and 
curvature. The equations for geodesics of a complex Finsler
space were first derived by H. Rund in \cite{Ru}; alternative presentations
are given in \cite{Pa}, \cite{Fa}, 
\cite{AP} and \cite{Sp}.  The equations for the Jacobi fields were  first 
determined by M. Abate and G. Patrizio in \cite{AP}.\par
The expressions given here have the peculiarity that, with no further 
arguing or manipulation,  
 they immediately reduce to the corresponding usual formulae 
of Hermitian geometry in case  $F$
is associated with an  Hermitian metric.\par
\medskip
 We believe that a careful
 investigation of the equations of Jacobi fields of the smoothly bounded convex domains in $\C^n$
can bring to isolate some crucial properties of the Kobayashi infinitesimal metric, 
which
characterize those domains up to biholomorphisms (for results in this direction, 
see e.g. \cite{AP}, \cite{BD}). 
A more detailed discussion of 
this topic will be the content of a forthcoming paper.
\bigskip
\bigskip
\subhead 2. Preliminaries and Notation
\endsubhead
\bigskip
In the whole paper, we will  use
 greek letters $\alpha$, $\beta$, etc.
for indices related to holomorphic vectors,  barred greek letters
$\bar\alpha$, $\bar \beta$, etc.  for indices related to the  
conjugated
vectors  and latin indices $i, j, k$, etc.
to denote  real vectors.\par
\medskip
$J_o$
is the  complex structure of $\C^n$ and 
 $<,>$  is the standard Hermitian  
product of $V = \C^n$. \par 
 The elements $\{\epsilon_0, \epsilon_1, \dots, 
\epsilon_{2n-1}\} \subset \C^n$ constitute
the standard real basis of $V= \R^{2n} = \C^n$ and they are  ordered 
so that 
$J_o(\epsilon_{2i}) = \epsilon_{2i + 1}$ for any $i = 0, \dots, n$.
We set
$\varepsilon_\alpha = \frac{1}{2}(\epsilon_{2\alpha} - \unim  
\epsilon_{2\alpha+1})$, $\alpha = 0, \dots, n-1$, and 
 $\varepsilon_{\bar \alpha} = 
\overline{\varepsilon_{\alpha}}$.
We also use the notation  $\{\epsilon^i\}$,
$\{\varepsilon^\alpha\}$ and $\{\varepsilon^{\bar\alpha}\}$
for the  dual bases of  $\{\epsilon_i\}$, $\{\varepsilon_\alpha\}$
and $\{\varepsilon_{\bar\alpha}\}$, respectively.\par
\medskip
We denote by $M$ a complex manifold with complex structure $J$. We also use the notation 
$\PTM = \TM/\C^*$, where $\TM = TM \setminus \{\text{zero section}\}$. \par
For any point $x\in M$
and any $v\in T_xM$, the tangent space $T_v(T_xM)$ is naturally identified with 
$T_xM$ and  we will use the  symbol $J$ also to denote the  complex structure 
on the tangent spaces $T_v(T_xM)$, given by the identification with $T_xM$. \par 
For any $v\in T_xM$, we denote by  $v^{10}$ and $v^{01}$
 the holomorphic and anti-ho\-lo\-mor\-phic parts w.r.t. $J$, that is: 
$$v^{10}= \frac{1}{2}(v - \sqrt{-1} J v)\ ,\qquad v^{01} = \overline{v^{10}} = 
\frac{1}{2}(v - \sqrt{-1} J v)\ .$$
\medskip
For any $x\in M$, a {\it linear frame \/} is an $\R$-linear 
isomorphism $u\: \R^{2n} \to T_xM$. A linear frame is called {\it complex linear frame\/}
 if it is a $\C$-linear isomorphism
$u\: \C^n = \R^{2n} \to T_xM$.
We  always 
identify a linear frame $u$  
with the corresponding basis $\{f_i\}$ in $T_xM$ defined by 
$$f_i = u(\epsilon_i)\in T_xM\ .$$
If a frame  $u$ is complex, we denote by $u^{10}$ the corresponding
 holomorphic basis, that is 
$$u^{10} = \{e_{\alpha}
 = u(\varepsilon_\alpha) = \frac{1}{2}(f_{2\alpha} - 
\unim f_{2\alpha+1})\}\ .$$ 
\medskip
If $(M, J, F)$ is a complex Finsler manifold, 
a complex linear frame $u = \{f_i\}$  is called {\it adapted\/}
if 
\roster
\item"a)" $f_0 \in S_x$ and $f_1 = J f_0$, where $S_x$
denotes the {\it Finsler pseudo-sphere\/} $S_x = \{\ v\in T_xM \ : F(v) = 1\ \}$; 
\item"b)" the vectors $f_2, \dots, f_{2n-1}$ span 
the maximal $J$-invariant
subspace  $\D_{f_0}$ of  $T_{f_0}S_x \subset T_{f_0}(T_xM) \simeq T_xM$;
\item"c)" the holomorphic vectors  $e_1$, \dots, $e_{n-1}$ constitute
 a unitary basis for $\D_{f_0}$ with respect to  
the Levi form $\L_x$ of $S_x$, corresponding to the defining function $\rho_F = F^2 - 1$. 
\endroster
The {\it unitary frame bundle\/} of $(M, J, F)$ is the subbundle $\UF \subset L^\C(M)$  given by all the 
adapted complex linear frames. \par
\medskip
It follows from definitions that any fiber of 
$\UF$ 
is invariant under the linear action of $U_n \times T^1 \subset GL_n(\C)$ on $T_xM$. 
Moreover, the orbit space of  this  action  can be  identified
with $\UF/U_n \times T^1 = \PTM$.\par
We will use the symbols  $\pi$, $\hat \pi$ and $\pi'$ to denote  the following
natural projections
$$\hat \pi: \UF \to \UF/U_{n-1}\times T^1 = \PTM \subset TM\ ,\qquad \pi': \PTM \to M\ ,$$
$$ \pi = \pi' \circ \hat \pi: \UF \to M\ .$$
\medskip
The {\it non-linear Hermitian connection of $\UF$\/} is  the unique 
distribution $\H$ on $\UF$,  which is complementary to the vertical distribution 
and which is invariant under the complex structure $\hat J$ of 
the complex linear frame bundle $L^\C(M)$. 
The distribution 
$\H$ is equal to the  real distribution underlying the 
holomorphic distribution of the real submanifold 
$\UF \subset L^\C(M)$ (see proof of Th. 3.9 in \cite{Sp}). \par
\medskip
The non-linear Hermitian connection $\H$ is uniquely 
determined by a connection form
$\omega$, that  is by a  $\gl$-valued 1-form on $\UF$ which verifies the following conditions:
\roster
\item"a)" for any $u\in \UF$, a vector $\X\in T_u\UF$ is so that $\omega_u(\X) = 0$ if and only if 
$\X \in \H_u$; 
\item"b)" if a vector $\X\in T_u\UF$ is vertical (i.e. $\pi_*(\X) = 0$), 
then $\omega_u(\X) = E_{\X}$, where 
$E_{\X}$ is the unique element in $\gl$ which generates an infinitesimal transformation 
on $L^\C(M)$ assuming the value $\X$ at the point $u$.\par
\endroster
Let 
$\omega^\alpha_\beta$, $\omega^{\bar \alpha}_{\bar \beta}$, $\theta^\alpha$ and $\theta^{\bar \alpha}$
be the components of the connection form and of the tautological 1-form in the  basis
$\{E^\alpha_\beta = \varepsilon_\beta\otimes \varepsilon^\alpha\}$
and $\{\varepsilon_\alpha\}$ of $\gl$ and of $\C^n$, respectively. In other words let  
$$\theta = \sum_{\alpha} \varepsilon_\alpha  \theta^\alpha
+ \sum_{\bar \alpha} \varepsilon _{\bar \alpha} \theta^{\bar \alpha}\ ,$$
$$\omega = \sum_{\alpha, \beta} E^\alpha_\beta \omega^\beta_{\alpha} + 
E^{\bar \alpha}_{\bar \beta} \omega^{\bar \beta}_{\bar \alpha}\ .$$
These 1-forms are not linearly independent but, 
at all points, they generate  the whole cotangent 
bundle $T^*\UF$. The
linear relations between them and the expressions for their
 exterior differentials are called {\it structure equations of the complex Finsler manifold 
$(M,J,F)$\/}.  We will shortly  list all such structure equations. \par
For this purpose 
 we first have to introduce some special $\C$-valued functions on $\UF$. 
\par
\medskip
For any  vector $X \in T_x M$,  denote by $V^X$
the vector field in $T(T_xM)$ which assumes the value $X$ at all points $U\in T_xM$.  
For any choice of vectors $X, Y, Z, W , U \in T_xM$, we define 
$$\bh_U(X,Y) = \left.V^X\left[V^Y\left(F^2\right)\right]\right|_U\ ;\qquad
\bH_U(X,Y,Z) = \left.V^X\left[V^Y\left[V^Z\left(F^2\right)\right]
\right]\right|_U\ ;\tag 2.1$$
$$
\gh_U(X,Y,Z,W) = \left.V^X\left[V^Y\left[V^Z\left[\hat 
V^W\left(F^2\right)\right]
\right]\right]\right|_U\ .
\tag2.2$$
For any adapted 
frame $u = \{f_i\}$ and  corresponding holomorphic frame
$u^{10} = \{e_\alpha\}$, we  set
$$h_{\alpha \beta}(u) = \bh_{f_0}(e_\alpha, e_\beta) \ ,\ 
H_{\alpha \beta \gamma}(u) = \bH_{f_0}(e_\alpha, e_\beta, e_\gamma)\ ,\  
H_{\alpha \beta \gamma\delta}(u) = \gh_{f_0}(e_\alpha, e_\beta, e_\gamma, e_\delta)
\ .$$
The symbols 
$h_{\bar \alpha \bar \beta}(u)$, $H_{\alpha \beta \bar\gamma}(u)$,
$H_{\alpha  \bar \beta \bar \gamma}(u)$, etc.
have analogous meanings. \par
Finally, in all following formulae, we will assume that 
the greek indices $\alpha$, $\beta$, $\gamma$, $\delta$, $\varepsilon$   
run between $0, \dots, n-1$; the indices
$\lambda, \mu, \nu, \rho, \sigma$ will instead run between $1$ and $n-1$.\par
\medskip
The first structure equations are given 
by the linear equations verified by the 
 1-forms $\omega^\alpha_\beta$ and $\omega^{\bar \alpha}_{\bar \beta}$: 
$$\omega^0_0 + \omega^{\bar 0}_{\bar 0} = 0\ ,\quad
\omega^0_\lambda +  \omega^{\bar \lambda}_{\bar 0} + h_{\lambda \nu} \omega^{\nu}_0
 = 0\ ,\quad
\omega^\lambda_\mu +
\omega^{\bar \mu}_{\bar \lambda} + 
H_{\bar \lambda\mu \nu} \omega^\nu_0 +
H_{\bar \lambda\mu\bar \nu} \omega^{\bar\nu}_{\bar 0} = 0\ .\tag2.3$$
In order to write down the 
expressions for the exterior differentials, it is 
convenient to replace the  1-forms $\omega^\alpha_\beta$, $\omega^{\bar \alpha}_{\bar \beta}$
with the following 1-forms
$\varpi^\alpha_\beta$ and $\varpi^{\bar \alpha}_{\bar \beta}$ 
$$\varpi^0_0 = \omega^0_0\ ,\ \ 
\varpi^\lambda_0 = \omega^\lambda_0\ ,\ \ 
\varpi^0_\lambda = - \omega^{\bar \lambda}_{\bar 0}\ , 
\ \ 
\varpi^\mu_\nu = \omega^\mu_\nu + 
H_{\bar \mu \nu \lambda} \omega^{\lambda}_0\ ,\ \ 
\varpi^{\bar \alpha}_{\bar \beta} = \overline{\varpi^\alpha_\beta}\ .\tag 2.4$$
Then the last structure equations are:
$$d\theta^\alpha + \varpi^\alpha_\beta \wedge \theta^\beta =
\Theta^\alpha + \Sigma^\alpha\ ;\tag2.5$$
$$d\varpi^0_0 + \varpi^0_\beta\wedge\varpi^\beta_0 = 
\Omega^0_0 \ ;\tag2.6$$
$$
d\varpi^\lambda_0 + 
\varpi^\lambda_\beta\wedge\varpi^\beta_0 =
\Omega^\lambda_0 + \Pi^\lambda_0\ 
,\quad
\ d\varpi^{0}_{\lambda} + 
\varpi^{0}_\beta\wedge\varpi^\beta_\lambda =
\Omega^0_\lambda + \Pi^0_\lambda\  ;
\tag2.7
$$
$$
d\varpi^\lambda_\mu + 
\varpi^\lambda_\beta\wedge\varpi^\beta_\mu = \Omega^\lambda_\mu
+ \Pi^\lambda_\mu + \Phi^\lambda_\mu\ ;
\tag2.8$$
where 
 $\Theta^\alpha$, $\Sigma^{\alpha}$, 
$\Omega^\alpha_\beta$, $\Pi^\lambda_0$,
$\Pi^0_\mu$, $\Pi^\lambda_\mu$ and 
$\Phi^\lambda_\mu$
are the following $\C$-valued 2-forms:
$$\Theta^\alpha = 
\frac{1}{2}
T^\alpha_{\beta \gamma} \theta^\beta\wedge\theta^\gamma\ ,\quad
\Sigma^\alpha = 
H_{\bar \alpha \mu \lambda}\varpi^\lambda_0\wedge \theta^\mu \ ,
\quad
\Omega^\alpha_\beta = 
R^{\alpha}_{\beta \gamma \bar \delta}\theta^\gamma\wedge\theta^{\bar \delta}
\ ,\tag 2.9$$
$$ \Pi^0_\lambda =  -\hat e_{\bar \gamma}(h_{\lambda \rho})\varpi^\rho_0 \wedge 
\theta^{\bar \gamma} \ ,
\ \ \Pi^\lambda_0 =  -\hat e_\gamma(h_{\bar \lambda \bar \rho})
\varpi^0_\rho \wedge 
\theta^{\gamma} 
\ , \tag2.10$$
$$
\Pi^\lambda_\mu= 
- \hat e_\gamma(H_{\bar \lambda \mu \bar \rho}) \varpi^{0}_\rho \wedge
\theta^\gamma - 
\hat e_{\bar \gamma}(H_{\bar \lambda \mu \rho})
\varpi^{\rho}_{0} \wedge\theta^{\bar \gamma}\ ,\tag2.11$$
$$\Phi^\lambda_\mu = 
\left( H_{\bar \lambda \bar \sigma \mu \rho }
-  h_{\bar \lambda \bar \sigma} h_{\mu \rho}  - 
 H_{\nu \bar \lambda  \bar \sigma} H_{\bar \nu \mu \rho }
 \right)
\varpi^\rho_0 \wedge \varpi^0_\sigma
\ .\tag2.12$$
for some suitable complex functions $T^\alpha_{\beta \gamma}$ and $R^{\alpha}_{\beta \gamma \bar \delta}$
on $\UF$.\par
\medskip
The 2-forms $\Theta$ and  $\Sigma$  are  called
{\it (pure) torsion form\/} and {\it 
Finsler torsion form\/}, respectively. 
The 2-form $\Omega$ is called the {\it (pure) curvature form\/}. 
A last,  we call $\Pi$ and $\Phi$  the
{\it oblique Finsler curvature\/} and  the {\it vertical
Finsler curvature\/}, respectively. \par
We recall that 
the Finsler curvature and torsion forms
are identically $0$ whenever the Finsler metric
is  associated with an Hermitian metric.
\bigskip
The next concepts will be essential for the  discussions
of the following sections, where also
the motivations for the terminology will appear clear.\par
\medskip
\definition{Definition 2.1} Let $(M,J, F)$ be a complex Finsler manifold 
and let $\varpi_H$ and $\varpi_K$ be the $\gl$-valued 1-forms on 
$\UF$ defined as 
$$\omega_H = \sum_{\lambda, \mu = 1}^n E^\mu_\lambda \omega^\lambda_{\mu} + 
E^{\bar \mu}_{\bar \lambda} \omega^{\bar \lambda}_{\bar \mu}\ ,
\quad 
\varpi_K = \sum_{\lambda, \mu = 1}^n E^\mu_\lambda \varpi^\lambda_{\mu} + 
E^{\bar \mu}_{\bar \lambda} \varpi^{\bar \lambda}_{\bar \mu}\tag 2.13$$
where $\varpi^\lambda_\mu$ and $\varpi^{\bar \lambda}_{\bar \mu}$
are the 1-forms defined in (2.4). Let also $\H$ the 
non-linear Hermitian connection on $U_F(M)$. Then the two distributions  on $\UF$
$\H'_H$ and $\H'_K$ 
defined by 
$$
\X \in \H'_H \ \Leftrightarrow\  \omega_H(\X) = 0\ ,\qquad
\X \in \H'_K \ \Leftrightarrow\  \varpi_K(\X) = 0\ ,$$
are called {\it semi-Hermitian connection\/} and 
{\it  Kobayashi connection for $(\PTM,J,F)$\/}, respectively.\par
Note that, at any $u\in \UF$, the subspaces $\H_H|_u$ and $\H_K|_u$
are both containing $\H_u$ as a proper subspace.
\enddefinition
\medskip
We conclude recalling the definition of connections and
Hermitian connections of complex vector bundles (see e.g. \cite{Ko1}).\par 
Let  $p : E \to N$
be a complex vector bundle over a manifold $N$ and 
denote by $\A^p$ denote the space of smooth $\C$-valued
p-forms on $N$. Denote also by $\A^p(E)$  the space of 
smooth complex  p-forms with values in $E$.
A {\it connection\/} $D$ on $E$ is a $\C$-linear homomorphism
$$D\ :\ \A^0(E) \longrightarrow \A^1(E)\ $$
such that 
$$D(f \sigma) = \sigma \otimes df + f \cdot D\sigma$$
for any $f\in \A^0$ and $\sigma \in \A^0(E)$.\par
In case $N$ is a complex manifold and $p: E\to N$ is a holomorphic vector bundle,  
a connection $D$ is called {\it holomorphic\/} if 
$$D^{01} = d^{01}$$
where $d$
is the usual exterior differential operator and 
$D^{01}:\A^0(E) \longrightarrow \A^{01}(E) $ and $d^{01}:
\A^0(E) \longrightarrow \A^{01}(E)$ are the 
components of $D$ and $d$, respectively, which transform
the sections $\sigma\in \A^0(E)$ into the 
 $(0,1)$-component  of 
$D(\sigma)$ and $d(\sigma)$. \par
In case $p: E \to N$ is an Hermitian vector bundle (i.e. endowed 
with a smooth family of Hermitian metrics on the fibers of $E$), a connection 
$D$ is called {\it Hermitian\/} if it is holomorphic and for any $\sigma, \rho \in 
\A^0$
$$d(g(\sigma, \rho)) = g(D\sigma, \rho) + g(\sigma, D\rho)$$
Recall that {\it on any Hermitian vector bundle there exists exactly 
one Hermitian connection\/}. 
\bigskip
\bigskip
\subhead 3. The  Hermitian  and  the
Kobayashi non-linear covariant derivatives of vector fields
\endsubhead
\bigskip
In this section, we introduce the definition of  covariant
derivation associated with the  distributions $\H_H$ and   $\H_K$ given in Definition 2.1. 
As mentioned in the Introduction, this covariant derivations can be defined using the
parallel transports along curves in $\PTM$, which are
determined by the curves in $\UF$ which are tangent to
the horizontal distribution $\H_H$ and $\H_K$, respectively (see also 
Rmk 3.8 in \cite{Sp}; note however that the discussion there concerns only 
curves $\gamma$ in $\PTM$ for which the vector $\pi'_*(\dot \gamma_t)) \in TM$
is nowhere vanishing). \par
However, we will adopt  here a different approach, which 
was considered by S. Kobayashi in \cite{Ko} and it is equivalent to the previous 
one, since it  is much more suitable 
for computations and further developments.\par
\medskip
Let us denote by $\TTM$ the  vector bundle $\TTM = (\pi')^{-1}(TM)$ defined as the pull-back bundle
w.r.t. the projection map $\pi'$. We thus obtain 
the following commuting diagram
$$\CD
\TTM @>{\tilde \pi'}>> TM\\
@V{\tilde \pi}VV  @VV{\pi}V\\
\PTM @>>{\pi'}> M
\endCD
\tag 3.1$$
It is clear that  there exists a unique complex structure on $\TTM$ (let us call it $\tilde J$), which makes 
$\tilde \pi: \TTM \to \PTM$ a holomorphic vector bundle. Moreover, 
as it was pointed out in \cite{Ko}, the Finsler metric $F$  on $(M, J)$
induces the following natural 
Hermitian metric on the vector bundle $\tilde \pi: \TTM \to \PTM$. \par
\medskip
Recall that the fiber $\tilde \pi^{-1}(v) \subset \TTM$ over an element $v\in \P T_xM$
coincides with the tangent space $T_xM$. Let $U\in T_xM$ be any non-zero vector 
which generates the 1-dimensional subspace $v = [U]$ and let $g_v$ be the bilinear form on 
$T_xM$ defined as 
$$g_v(X,Y) = 2 \operatorname{Re}(\bh_U(X^{10}, Y^{01})) = 
\frac{1}{2} \left(
h_U(X,Y) + h_U(JX, JY)\right)\tag3.2$$
where $\bh_U$ is as in (2.1). It is clear that right hand side of (3.2)
 is an Hermitian metric on $T_xM$. Moreover
from the invariance 
properties of complex Finsler metrics under $\C^*$-multiplications, 
it follows that for any $\lambda \in \C^*$
$$h_{\lambda V}(X,Y) + h_{\lambda V}(JX, JY)  = h_{V}(X,Y) + h_{V}(JX, JY)\tag 3.3$$
(see e.g. \cite{Sp} Lemma 2.4 b) and c), or \cite{Ko}). This
 shows that the r.h.s. of (3.2) is indeed an Hermitian metric which depends only on 
the line  $v = [U] \in \P T_xM$.\par
\medskip
 Let us now denote by $\CTTM$ the set of all local sections
$X: \PTM \to \TTM$. Notice that  for any $X\in \CTTM$ 
there exists 
at least one local vector field $\X$ on $\UF$, such that 
$$\pi_*(\X_u) = X(v)\ ,$$
for any $v\in \P T_xM$
and any $u\in \hat \pi^{-1}(v) \subset \UF$. If this is the case, 
we will say that {\it $\X$ is $TM$-projectable\/} (or, more often, just {\it projectable\/})
 and we will call $X$ the  
{\it projection of $\X$ in $\CTTM$\/}.\par
Similarly, for any local vector field $\hat X$ on $\PTM$, there exists 
some local vector fields $\X$ on $\UF$, such that for any frame $u\in \UF$, 
$$\hat \pi_*(\X_u) = \hat X_{\hat \pi(v)}\ .$$
In this case we will say that {\it $\X$ is $\PTM$-projectable\/} 
(or  just {\it projectable\/}) and  we will call $\hat X$ the
{\it the projection of  $\X$ in $\PTM$\/}.\par
\medskip
Now,  the following  technical lemma is required. \par
\bigskip
\proclaim{Lemma 3.1} 
Let  $\hat X$ and $Y$ be  a (local) vector field of $\PTM$ and 
a (local) section in $\CTTM$, respectively, and let  
 $\X$ and $\Y$ two projectable  vector fields on $\UF$, such that 
$X$ is the projection  of $\X$ on $\PTM$ and  $Y$ is the projection of $\Y$ in $\CTTM$.\par
The functions  on $\UF$
$$F^{\X, \Y}(u) = u\left(\underset{\phantom{A}}\to{\overset{\phantom{A}}\to\X_u(\theta(\Y))}
 + \omega_u (\X) \cdot \theta_u(\Y)\right)\tag3.4$$
$$G^{\X, \Y}(u) = u\left(
\underset{\phantom{A}}\to{\overset{\phantom{A}}\to
\X_u(\theta(\Y))
} + \varpi_u (\X) \cdot \theta_u(\Y)\right)\tag3.5$$
assume constant values along the fibers $\hat \pi^{-1}(v) \in \UF$
and are independent  on the choice of the projectable vector fields $\X$ and $\Y$.\par
In particular, $F^{\X, \Y}$ and $G^{\X, \Y}$ define elements of
$\CTTM$, which  depend linearly on the value $\hat X_v$, at any  $v\in \PTM$.
\endproclaim
\demo{Proof} In order to prove that 
$F^{\X, \Y} (u)$ is constant along $\pi^{-1}(v)$,  it suffices to 
check that $\tilde A_u\left(F^{\X, \Y}\right) \equiv 0$ for any $u\in \pi^{-1}(v)$ and 
any vertical vector field $\tilde A$ on $\UF$. Indeed, we may   consider 
only vertical vector fields
$\tilde A$ which are fundamental vector fields, associated with elements $A \in \u$
(for the Def.  of fundamental vector fields, see e.g. \cite{KN} vol.I). \par
Notice that
$$\tilde A\left(\X(\theta(\Y)) + \omega(\X) \cdot \theta(\Y)\right) = $$
$$ =
[\tilde A, \X](\theta(\Y)) + \X\left( \L_{\tilde A} \theta(\Y)\right) + 
\X\left(\theta([\tilde A, \Y])\right) + 
\L_{\tilde A}\omega(\X) \cdot \theta(\Y) +$$
$$ +
\omega([\tilde A, \X]) \cdot \theta(\Y) + 
\omega(\X) \cdot \left(\L_{\tilde A} \theta(\Y)\right) + 
\omega(\X) \cdot \theta([\tilde A, \Y])\ . \tag3.6$$
On the other hand 
$$\hat \pi_*([\tilde A, \Y]) = [\hat \pi_*(\tilde A), \hat \pi_*(\Y)] = 0\qquad 
\hat \pi_*([\tilde A, \X]) = [\hat \pi_*(\tilde A), \hat \pi_*(\X)] = 0\ .$$
 In particular, 
$[\tilde A, \Y]$
and $[\tilde A, \X]$ are both vertical vector fields for the bundle 
$\hat \pi: \UF \to \PTM$ and we can write that
$[\tilde A, \X]_u = \tilde B_u$,  where $\tilde B$  is a fundamental vector
field associated with  an element  
 $B\in \u$. This implies that
$$\theta([\tilde A, \Y]) = 0\ ,\quad [\tilde A, \X](\theta(\Y))_u = 
- B \cdot \theta(\Y)_u =  - \omega([\tilde A, \X]) \cdot \theta(\Y)_u\ .$$
Then  (3.6) becomes equal to 
$$\tilde A(\X(\theta(\Y) + \omega(\X) \cdot \theta(\Y)) = $$
$$ = - \omega([\tilde A, \X]) \cdot \theta(\Y)_u - 
A \cdot \X(\theta(\Y)) - A\cdot \omega(\X) \cdot \theta(Y) + 
\omega(\X) \cdot A \theta(Y) +$$
$$ + \omega([\tilde A, \X]) \cdot \theta(\Y) - 
\omega(\X) \cdot A\cdot \theta(\Y) = 
 - A\cdot \left[\X(\theta(\Y)) + \omega(\X) \cdot \theta(Y)\right]\ .\tag 3.7$$
From (3.7) and the definition of $F^{\X, \Y}$, 
it follows immediately that $\tilde A_u\left(F^{\X, \Y}\right) \equiv 0$
for any $u\in \pi^{-1}(v)$.\par
\medskip
Now, consider other two projectable  vector fields $\X'$, $\Y'$, of which 
$X$ and $Y$ are the corresponding projections. Then, using the structure 
equations (2.5) - (2.8) we have
$$F^{\X', \Y'} (u) - F^{\X, \Y} (u) = $$
$$ =
(\X' - \X)(\theta(\Y')) + \omega(\X' - \X) \cdot \theta(\Y') + 
\X(\theta(\Y' - \Y)) + \omega(\X) \cdot \theta(\Y' - \Y) = $$
$$ = (\X' - \X)(\theta(\Y')) + \omega(\X' - \X) \cdot \theta(\Y') = $$
$$ = d\theta(\X' - \X, \Y') + \omega(\X' - \X) \cdot \theta(\Y') = $$
$$ = - \varpi(\X' - \X) \cdot \theta(\Y')  + 
 - \varpi(\Y') \cdot \theta(\X' - \X) + 
\Theta(\X' - \X, \Y') + $$
$$ +\Sigma(\X' - \X, \Y) + \omega(\X' - \X) \cdot \theta(\Y')$$
At this point, we remark  that $\varpi(\X' - \X) = \omega(\X' - \X)$: in fact
$\hat \pi_*(\X' - \X) = 0$ and this implies that the 1-forms $\omega^\lambda_0$ and 
$\omega^0_{\lambda}$ vanish on $(\X' - \X)$, by (5.32) in \cite{Sp}. 
The same argument implies that
$\Theta(\X' - \X, Y) = \Sigma(\X' - \X, \Y) = 0$ and hence
$$F^{\X', \Y'}(u) - F^{\X, \Y} (u) = 0\ .$$
This concludes the proof of both claims for the function $F^{\X, \Y}$. 
The proof of the corresponding claims for the function $G^{\X, \Y}$ is based
on very similar arguments.\qed
\enddemo
\bigskip
By means of Lemma 3.1, the following objects are well defined.\par
\bigskip
\definition{Definition 3.2} For any 
$Y\in \CTTM$, let $\nabla Y$ and $D Y$ be the 
elements in $\A^1(\TTM)$ defined by 
$$\derI {\hat X} Y v = 
u\left(\underset{\phantom{A}}\to{\overset{\phantom{A}}\to
\X(\theta(\Y))}+ \omega_u(\X)\cdot \theta_u(\Y)
\right)\ ,\tag3.8$$
$$\derII {\hat X} Y v= u\left(\underset{\phantom{A}}\to{\overset{\phantom{A}}\to
\X(\theta(\Y))} + \varpi_u(\X)\cdot \theta_u(\Y)
\right)\ ,\tag3.9$$
for any vector field $\hat X$ on $\PTM$ and any $v \in \PTM$; here 
$u$ is any frame  of $\hat \pi^{-1}(v) \subset \UF$ and $\X$, $\Y$
are two projectable  vector fields of $\UF$, whose projections are 
$\hat X$ and $\Y$, respectively. \par
 $\nabla$ and $D$ are connections on the holomorphic vector bundle
$\TTM$. For any section $Y \in \CTTM$, we call $\derI {\hat X} Y v$
and $\derII {\hat X} Y v$  the {\it non-linear 
semi-Hermitian covariant derivative\/} 
and the {\it non-linear Kobayashi covariant derivative\/}, respectively, 
along $\hat X$ at the point $v$.
\enddefinition
\bigskip
\remark{Remark 3.3} The set of local 
vector fields on $M$ can be naturally identified with the element in 
$\CTTM$, which are  sections that are constant along the 
fibers of  $\pi': \PTM \to M$. \par
 From this  and the properties of the distribution $\H$, 
it follows that 
whenever $F$ is associated with 
an Hermitian metric and $Y$ is a vector field on $M$, then 
both non-linear covariant derivatives $\derI {\hat X} Y v$
and $\derII {\hat X} Y v$ depend only on 
 $X = \pi'_*(\hat X) \in TM$ and on  $x = \pi'_*(v) \in M$, and
they both coincide with  usual {\it linear} Hermitian covariant derivative
of  $Y$ along  $X$.
\endremark
\bigskip
In the following Proposition, we give two  useful characterizations of 
the connections $\nabla$ and $D$. In particular  we show that $D$
coincides with the Finslerian connection introduced by S. Kobayashi in \cite{Ko}.\par 
\bigskip
\proclaim{Proposition 3.4} Let $(M, J, F)$ be a complex Finsler manifold, 
$\nabla$ and $D$  as in Definition 3.2 and $g$ the Hermitian metric on 
$\TTM$ defined in (3.2). \par
 Then for any local
vector field $\hat X$ on $\PTM$ and  any  sections
$Y, Z\in \CTTM$  
$$\derI {\hat X} {JY} v = J\derI {\hat X} {Y} v\ ,\qquad 
\derII {\hat X} {JY} v = J\derII {\hat X} {Y} v\ ,\tag3.10$$
$$\hat X (g(Y, Z))|_v + g_v(\derI {\hat X} {Y} v, {Z}) +
g_v(Y, \derI {\hat X} {Z} v) =$$
$$ = \bH_U(X,Y^{10},Z^{01}) + \bH_U(X,Z^{10},X^{01})\ ,\tag3.11$$
$$\hat X (g(Y, Z))_v + g_v(\derII {\hat X} {Y} v, {Z}) +
g_v(Y, \derII {\hat X} {Z} v) = 0\ .\tag3.12$$
where $U$ is any non-trivial vector in the complex line $v = [U]\in \P T_xM$, 
$x = \pi'(v)$ and $H$ is the trilinear function 
defined in (2.1).\par
Furthermore, the connection $D$ is  holomorphic  and 
 it coincides with the Hermitian connection 
of the  Hermitian bundle $\tilde \pi: \TTM \to \PTM$.
\endproclaim
\demo{Proof}
The Hermitian metric $g_v$ can be conveniently expressed using 
the components $\theta^\alpha$ of the tautological 1-form of $\UF$. In fact, 
for given two local sections $X, Y$ in $\CTTM$, consider
 two projectable vector fields $\X$ and $\Y$ on $\UF$, which project onto 
$X$ and $Y$; then for any $v\in \PTM$,
$$g_v(X,Y) = \sum_{\alpha} (\theta^\alpha(\X) \theta^{\bar \alpha}(\Y) +
 \theta^\alpha(\Y) \theta^{\bar \alpha}(\X))_u = <\theta_u(\X), \theta_u(\Y)>\tag 3.13$$
where $u$ is any frame in $\hat \pi^{-1}(v)$.\par
\medskip
 Consider now three vector fields $\X$, $\Y$ and $\cZ$ on $\UF$ which 
project onto $\hat X$, $Y$ and $Z$, respectively. Let also
$\Y^J$ be a vector field on $\UF$ which projects onto the local section $J Y \in \CTTM$. Then
$$\derI {\hat X} {JY} v = u(\X(\theta(\Y^J)) + \omega(\X)\cdot \theta(\Y^J)) = 
u(J_o \X \theta(\Y) + \omega(\X) \cdot [J_o \theta(\Y)] ) = $$
$$ = J  u(\X(\theta(\Y^J)) + \omega(\X)\cdot \theta(\Y))  = 
J \derI {\hat X} Y v\ .$$
In a similar way one can prove that $\derII {\hat X} {JY} v = J \derII {\hat X} {Y} v$.\par
Now, for  (3.11), one should observe that
$$\hat X (g(Y, Z))_v = \X(<\theta_u(\Y), \theta_u(\cZ)>) = $$
$$ =
<\X(\theta_u(\Y)), \theta_u(\cZ)> + <\theta_u(\Y), \X\theta_u(\cZ)> = $$
$$ = g_v(\derI {\hat X} Y v, Z) + g_v(Y, \derI {\hat X} Z v) - $$
$$ -
<\omega_u(\X)\cdot \theta_u(\Y), \theta_u(\cZ)> - 
<\theta_u(\Y), \omega_u(\X)\cdot \theta_u(\cZ) > $$
On the other hand,
$$<\omega_u(\X)\cdot \theta_u(\Y), \theta_u(\cZ)> +  <\theta_u(\Y), \omega_u(\X)\cdot \theta_u(\cZ) > = $$
$$ =
\sum_{\alpha, \beta} \omega^\alpha_\beta(\X)\theta^\beta(\Y)\theta^{\bar \alpha}(\cZ) + 
\overline{\omega^\beta_\alpha(\X)}\theta^\beta(\Y)\theta^{\bar \alpha}(\cZ) + $$
$$ +
 \omega^\alpha_\beta(\X)\theta^\beta(\cZ)\theta^{\bar \alpha}(\X) + 
\overline{\omega^\beta_\alpha(\X)}\theta^\beta(\cZ)\theta^{\bar \alpha}(\X)$$
By (2.3), we get 
$$<\omega_u(\X)\cdot \theta_u(\Y), \theta_u(\cZ)> +  <\theta_u(\Y), \omega_u(\X)\cdot \theta_u(\cZ) > = $$
$$ - \bH_U(\X,Y^{10}, Z^{01}) - \bH_U(\X,Y^{01}, Z^{10})\ ,$$
which proves (3.11).
(3.12) can be proved in the same way, using 
the equations (2.4) in place of (2.3). \par
\medskip
To conclude, we have to show that for any $Y \in \CTTM$ and any vector field $\hat X$ in $\PTM$ 
$$\derII {\hat X + \unim J \hat X} Y {} = (\hat X + \unim J \hat X)(Y)\ .$$ 
One can verify that this condition is equivalent to show that for 
any vector field $\hat X$ in $\PTM$  there exist two projectable vector 
fields  $\X$, $\X_J$ on $\UF$ which project onto $\hat X$ and $J \hat X$, respectively, 
and so that, for any $\lambda, \mu  = 1, \dots , n-1$,
$$\varpi^\lambda_\mu(\X) =  \varpi^0_0(\X) = 0 
\ ,\qquad \varpi^\lambda_\mu(\X_J) =  \varpi^0_0(\X_J) = 0\ .\tag 3.14$$
From (2.4), it is clear that the distribution $\Cal B$ defined 
by the conditions (3.14) consists of the set of  vector spaces 
$$\Cal P_u = \H_u \oplus \tilde{\Cal V}_u \subset T_u\UF\ ,$$
where  $\H$ is the non-linear Hermitian connection of $\UF$ and 
$$\tilde{\Cal V}_u \= \{ \ \X_u\ :\ \pi_*(\X_u) = 0\ ,
\ \omega^\lambda_\mu(\X_u) = - H_{\bar \lambda \mu \nu} \omega^\nu_0(\X_u)\}\ .$$
By the results in \cite{Sp}, one can check that the distribution $\Cal B$ is invariant under the action 
of $U_{n-1} \times T^1$ and, at all points, there exists a complex structure 
$\tilde  J_u : \Cal P_u \to \Cal P_u$ which is $U_{n-1} \times T^1$-invariant and 
 projects onto the complex structure of $T_{\hat \pi(u)}(\PTM)$. In fact, 
$\Cal B$ is spanned by the vector fields 
 $\re(\hat e_\alpha)$, $\im(\hat e_\beta)$, $\re(\tilde e'_\lambda)$, 
$\im(\tilde e'_\mu)$, defined in \S 5.1 - 5.2 in \cite{Sp}. Therefore, 
by Prop. 5.4 and Prop. 5.5 (1)
in \cite{Sp}), it follows that $\Cal B$ and the 
(almost) complex structure on $\Cal B$ defined by
$$\tilde J_u(\re(\hat e_\alpha)) = \im(\hat e_\beta)\ ,\qquad 
\tilde J_u(\re(\tilde e_\lambda)) = \im(\hat e_\lambda)\ , $$
are both invariant  under the action of $U_{n-1} \times T^1$.
\par
To conclude the proof, it is enough to take as vector fields $\X$ and $\X_J$
there  unique 
vector field  $\X$ on $\Cal B$, which projects
onto $\hat X$,  and that the vector field $\X_J = \tilde J \X$, respectively. \qed
 \enddemo
\bigskip
\bigskip
\subhead 4. Torsion and  curvature of the Kobayashi connection
\endsubhead
\bigskip
We now want  to define the torsion and the curvature of the Kobayashi connection 
of the vector bundle $\tilde \pi: \TTM \to \PTM$ and express them 
in terms of the Finsler torsions and 
curvatures of the non-linear Hermitian connection on $\UF$. This can be 
done by virtue of the following proposition.\par
\bigskip
\proclaim{Proposition 4.1} Let $\hat X, \hat Y$  be  local
vector fields on $\PTM$ and  $Z\in$ $\CTTM$. Let also 
$X, Y$ be the sections in $\CTTM$ defined by 
$$X(v) = \pi'_*(\hat X(v))\ ,\qquad  Y(v) = \pi'_*(\hat Y(v))\ .$$
For any $v\in \PTM$, consider the vectors in $T_{\pi'(v)} M$ defined by
$$T_{\hat X, \hat Y}(v) = \derII{\hat X} Y v - \derII{\hat Y} X v - [X,Y]\ .\tag4.1$$
$$R_{\hat X, \hat Y} Z (v) = \derII{\hat X}{(\derII{\hat Y} Z {\tilde v})} v - 
\derII{\hat Y}{(\derII{\hat X} Z {\tilde v})} v - \derII{[\hat X, \hat Y]} Z v\tag 4.2$$
Then if  
$\X$, $\Y$ and $\cZ$  are  three projectable vector fields
on $U_F(M)$, which project onto  $\hat X$,  
$\hat Y$ and $Z$, respectively,  and if 
$u: \PTM \to \UF$ is any local section of the bundle $\hat \pi: \UF \to \PTM$, then 
$T_{\hat X, \hat Y}(v)$ and $R_{\hat X, \hat Y} Z (v)$ verify
$$T_{\hat X, \hat Y}(v)  = u_v(\Theta(\X,\Y) + \Sigma(\X,\Y))\tag4.3$$ 
$$R_{\hat X, \hat Y} Z (v) = u_v(\Omega(\X, \Y)\cdot \theta(\cZ) + \Pi(\X, \Y) \cdot \theta(\cZ) + 
\Phi(\X, \Y)\cdot \theta(\cZ))\ .\tag 4.4$$
In particular,  $T_{\hat X, \hat Y}(v)$ and $R_{\hat X, \hat Y} Z (v)$
depend only on the values  $\hat X_v, \hat Y_v$ and $Z_v$.
\endproclaim
\demo{Proof} By definitions
$$u_v^{-1}(\derII{\hat X}{Y} {v} - 
\derII{\hat Y}{X} {{v}} -  \pi'_*([\hat X,\hat Y]_v)) = $$
$$\X(\theta_{u}(\Y))|_{u_v} - \Y(\theta_u(\X))|_{u_v} +
\varpi_{u_v}(\X) \cdot \theta_{u_v}(\Y) - \varpi_{u_v}(\Y) \cdot \theta_{u_v}(\X) - 
\theta([\X, \Y])_{u_v} = $$
$$ = d\theta_{u_v}(\X, \Y) -  \varpi_{u_v}(\X) \cdot \theta_{u_v}(\Y) - \varpi_{u_v}(\Y) 
\cdot \theta_{u_v}(\X) = $$
$$ = \Theta_{u_v}(\X, \Y) + \Sigma_{u_v}(\X, \Y)\ ,$$
and this proves (4.3). Similarly
$$u^{-1}_v(\derII{\hat X}{(\derII{\hat Y} Z v)} {v} - 
\derII{\hat Y}{(\derII{\hat X} Z v)} {v} - \derII{[\hat X, \hat Y]} Z {v}) = $$
$$ = [\X(\varpi(\Y)) - \Y(\varpi(\X)) - \varpi([\X,\Y])]\cdot \theta(\cZ) + 
+ \varpi(\Y)\cdot \X(\theta(\cZ)) - \varpi(\X)\cdot \Y(\theta(\cZ)) + $$
$$ + \varpi(\X) \cdot \Y(\theta(\cZ)) - \varpi(\Y) \cdot \X(\theta(\cZ)) +
\varpi(\X) \cdot \varpi(\Y) \cdot \theta(\cZ) - \varpi(\Y) \cdot \varpi(\X) \cdot \theta(\cZ) = $$
$$ = d\varpi(\X, \Y) \cdot \theta(\X, \cZ)  +
\varpi(\X) \cdot \varpi(\Y) \cdot \theta(\cZ) - \varpi(\Y) \cdot \varpi(\X) \cdot \theta(\cZ) \ . $$
Now, using the structure equations (2.7) and (2.8), it follows that
$$d\varpi(\X, \Y) \cdot \theta(\X, \cZ)  +
\varpi(\X) \cdot \varpi(\Y) \cdot \theta(\cZ) - \varpi(\Y) \cdot \varpi(\X) \cdot \theta(\cZ) = $$
$$ = [\Omega_{u_v}(\X, \Y) + \Pi_{u_v}(\X, \Y) + \Phi_{u_v}(\X, \Y)]_u\cdot \theta_{u_v}(\cZ) = 
R(X,Y Z, u_v)\ ,$$
and this proves (4.4).\qed
\enddemo
\bigskip
By means of Proposition 4.1, we may define the torsion and  the curvature  of the Kobayashi connection
as follows
(see also \cite{KN} or
\cite{Ko1}). 
\medskip
\definition{Definition 4.2} The {\it torsion of
Kobayashi connection on $\TTM$\/}  is the element in $\Cal A^2(\TTM)$, defined by 
$$T_{\hat X, \hat Y}(v) = \derII{\hat X} Y v - \derII{\hat Y} X v - [X,Y]\tag4.5$$
for any  vector fields $\hat X, \hat Y$ on $\PTM$; 
here $X$ and $Y$ are the sections in $\CTTM$ defined by 
$X(v) = \pi'_*(\hat X(v))$, $Y(v) = \pi'_*(\hat Y(v))$.\par
The {\it curvature of the Kobayashi connection on $\TTM$ \/} is the $\C$-linear operator
$$R\ :\ \CTTM \longrightarrow \Cal A^2(\TTM)$$ 
defined by 
$$R_{\hat X, \hat Y} Z (v) = \derII{\hat X}{(\derII{\hat Y} Z {\tilde v})} v - 
\derII{\hat Y}{(\derII{\hat X} Z {\tilde v})} v - \derII{[\hat X, \hat Y]} Z v\tag 4.6$$
for any vector fields $\hat X, \hat Y$ on $\PTM$  and any $Z\in \CTTM$
\enddefinition
\medskip
Note that, by Proposition 4.2, 
the values $T_{\hat X, \hat Y}(v)$ and $R_{\hat X, \hat Y} Z (v)$ depends only 
on the values $\hat X|_v$, $\hat Y|_v$ and $Z|_v$.\par
\bigskip
\bigskip
\subhead 5. Equations of geodesics
\endsubhead
\bigskip
We are going to write down the  equations in terms of the Kobayashi connection. 
It will be a simple corollary of the results of \cite{Sp} and of the previous discussion.\par
\medskip
Let $\gamma: [a, b] \to M$ be 
a regular curve in $M$ and  $v_\gamma$ the corresponding
curve in $\PTM$ defined by 
$$v_\gamma: [a, b] \longrightarrow \PTM\ ,$$
$$v_\gamma(t) = [\dot \gamma_t]\in\Bbb P T_{\gamma_t}M$$
It is not difficult to realize that any
tangent  vector $\doubledot \gamma(t)$ 
depends linearly on the  second derivative $\ddot \gamma_t$
and non-linearly on the first derivative $\dot \gamma_t$. \par
To simplify the notation, in the following 
 we will  assume that  the three 
functions $\dot \gamma$, $v_\gamma$ and $\doubledot \gamma$ are always evaluated at 
the same point $t\in [a,b]$.\par
\medskip
 We say that
 a complex Finsler manifold $(M, J,F)$ is 
 {\it geodetically torsion free\/} if for any 
$v\in \PTM$ and any $0\neq \hat X, \hat U\in T_v\PTM$, with $\pi_*(\hat U) = U$ so that
 $U$ is 
a non trivial vector in $v = [U]$,
$$g_{v_\gamma}( T_{\hat X \hat U}, 
U) = 0\ .\tag 5.1$$
Note that from (3.10), (4.3)  and the definitions of Finsler torsions $\Theta$ and $\Sigma$, 
a complex Finsler manifold is geodetically torsion free if and 
only if it is geodetically torsion free w.r.t.   Def. 6.5 of \cite{Sp}.\par
The geodetically torsion free complex Finsler manifolds coincide with the 
manifold called  {\it weakly K\"ahler\/} 
by Abate and Patrizio in \cite{AP}. This term is motivated by the fact that, 
whenever a Finsler metric 
is associated with an Hermitian metric, it is geodetically torsion free
if and only if the corresponding Hermitian metric is torsion free and hence K\"ahler.\par
\medskip
\proclaim{Theorem 5.1} 
Let $(M,J,F)$ a  complex 
Finsler manifold. 
Then  $\gamma$ is a geodesic of $M$ if and only if 
$$\derII{{\doubledot {\gamma}}}{\dot \gamma} {v_{\gamma}} + g_{v_\gamma}( T_{\doubledot \gamma V_\gamma}, 
V_\gamma) \equiv 0\ .\tag5.2$$
at any $t \in [a,b]$, where $V_\gamma$ is any  vector in $T\P T_{v_\gamma}M$
such that $[\pi'_*(V_\gamma)] = v_\gamma$. In particular, if 
$(M, J, F)$ is geodetically torsion free, then $\gamma$ is a geodesic if and only if 
$$\derII{{\doubledot {\gamma}}}{\dot \gamma} {v_{\gamma}} = 0\ .\tag5.3$$
\endproclaim
\demo{Proof}  From [Sp]  Th. 6.2, 
we have that $\gamma$ is a geodesic if and only if 
for any lift $\tilde \gamma : [a,b] \to \UF$ 
$$\dot{\tilde \gamma_t}(\theta^0(\dot{\tilde \gamma_t})) + 
\varpi^0_0(\dot{\tilde \gamma_t}) \theta^0(\dot{\tilde \gamma_t}) = 0\ ,\tag 5.4$$
$$\dot{\tilde \gamma_t}(\theta^{\bar 0}(\dot{\tilde \gamma_t})) + 
\varpi^{\bar 0}_{\bar 0}(\dot{\tilde \gamma_t})  \theta^{\bar 0}(\dot{\tilde \gamma_t}) = 0\ ,
\tag 5.5$$
$$\varpi^0_\lambda(\dot{\tilde \gamma_t}) + T^0_{\lambda 0} \theta^0(\dot{\tilde \gamma_t})
 = 0 = \varpi^{\bar 0}_{\bar \lambda}(\dot{\tilde \gamma_t})  + \overline{T^0_{\lambda 0} \theta^0
(\dot{\tilde \gamma_t})}\ .
\tag 5.6$$
We recall that a curve $\tilde \gamma: [a,b] \to \UF$ is a lift of $\gamma$
if and only if $\pi\circ \tilde \gamma = \gamma$ and for any 
frame $\tilde \gamma_t = \{f_0(t), \dots , f_{2n-1}(t)\}$, the vector 
$f_0(t)$ belongs to $f_0(t) \in \C^* \dot \gamma_t$. \par
From the definition, it follows that $\theta^\lambda(\dot{\tilde \gamma_t}) = 
\theta^{\bar \lambda}(\dot{\tilde \gamma_t}) = 0$ for any $t$ and hence
(5.4) - (5.6) are equivalent to
$$\dot{\tilde \gamma_t}(\theta(\dot{\tilde \gamma_t}))+ 
\varpi(\dot{\tilde \gamma_t}) \cdot \theta(\dot{\tilde \gamma_t}) + \Theta^0(
\dot{\tilde \gamma_t}, e_0) = 0\tag5.7$$
and this implies (5.7).\par
Conversely, if (5.1) holds for any $t$, it follows immediately that (5.7)
holds for any lift and hence $\gamma$ is a geodesic, by \cite{Sp} Th. 6.2.\qed
\enddemo
\bigskip
\bigskip
\subhead 6. Jacobi fields
\endsubhead
\bigskip
Let $\gamma: [a, b] \to M$ be a geodesic of 
$(M, J, F)$ and let $V : (- \delta, \delta) \times [a,b] \to M$ be  a smooth map 
such that  $V(s, *) = \gamma^{(s)}: [a, b] \to M$ is a geodesic for any 
$s\in [-\delta,\delta]$ and with $\gamma^{(0)} = \gamma$. We call $V$ a {\it 1-parameter
family of geodesics centered at $\gamma$\/}. We recall that a
vector field $I$ on $\gamma([a,b])$ is called {\it Jacobi field
for $\gamma$\/} if and only if it is of the form
$$I_{\gamma_t} = \left.\frac{d}{ds}(V(*, t))\right|_{s=0}$$
for some 1-parameter family of geodesics centered at $\gamma$. \par
The goal of this section is to determine the differential equations
which characterize the Jacobi fields, using the Kobayashi connection and its
torsion and curvature.\par
\medskip
Let $X$ be a  vector field defined on the points 
of a curve $\gamma: [a,b] \to M$. We call {\it standard lift of $X$ along $\gamma$
\/} the vector field  $L_{(X,\gamma)}$ along the curve $v_\gamma: [a,b] \to \PTM$
defined as follows. \par
Extend $X$ to a local vector field and let
$\Phi_s^{X}: M \to M$ the corresponding flow. Then let 
$$V: (-\delta, \delta) \times [a,b] \to M\ ,
\qquad V(s, t)  = \Phi^X_s(\gamma_t)$$ 
and let $v_{\gamma^{(s)}}(t)$ the 1-parameter family of lifted curves
$$v_{\gamma^{(s)}}(t) = [\dot \gamma^{(s)}_t]\in \PTM$$
where  $\gamma^{(s)}_t = V(s,t)$. We set
$$L_{(X,\gamma)}|_{v_\gamma} \= \left.\frac{d}{ds} v_{\gamma^{(s)}}(t)\right|_{s=0}\ .$$
\medskip
\proclaim{Theorem 6.1} Let $I$ be a vector field defined on the points
of a geodesic $\gamma$. Then $I$ is a Jacobi vector field if and only 
if it verifies the following system of
equations at all points of the geodesic:
$$\derII {\doubledot \gamma} {\left(\derII {\doubledot \gamma} I {v_{\gamma}}\right )}{v_\gamma} -
 R_{\doubledot \gamma L_{(I,\gamma)}} \dot \gamma  -  \derII {\doubledot \gamma} {\left(T_{
\doubledot\gamma L_{(I,\gamma)}}(v_\gamma)\right)}
{v_\gamma} = 0\ .\tag6.1$$
\endproclaim
\demo{Proof} Let $V : (- \delta, \delta) \times [a,b] \to M$ be  a 1-parameter
family of geodesics centered at $\gamma$  so that
$$I_{\gamma_t} = \left.\frac{d}{ds}(V(*, t))\right|_{s=0}\ .$$
Let also $\hat V: (- \delta, \delta) \times [a,b] \to \PTM$ be the associated
map such that 
$$\hat V(s, t) = [v_{\gamma^{(s)}}(t)]\ .$$
Then we may consider the vector fields $\hat X = \hat V_*(\frac{\partial}{\partial t})$
and $\hat Y = \hat V_*(\frac{\partial}{\partial s})$ and the associated functions on 
$(- \delta, \delta) \times [a,b]$ with values in $TM$ defined by 
$$X(s,t) = \pi'_*(\hat X(s,t))\ ,\qquad Y = \pi'_*(\hat Y(s,t))\ .$$
 Clearly, $[\hat X, \hat Y] = 0$
as well as $[X,Y] = 0$. Moreover,
$$\hat X_{v_\gamma(t)} =\doubledot \gamma(t)\ ,\qquad \hat Y_{v_{\gamma}(t)} = L_{(I,\gamma)}\ ,
\qquad X_{\gamma(t)} = \dot \gamma_t\ ,\qquad Y_{\gamma(t)} = I_{\gamma_t}\ .$$
Therefore
$$\derII {\doubledot \gamma} {\left(\derII {\doubledot \gamma} I {v_{\gamma}}\right )}{v_\gamma} = 
\derII {\hat X} {\left(\derII {\hat X} Y {v_{\gamma}}\right )}{v_\gamma} = $$
$$ \derII {\hat X} {\left(\derII {\hat Y} X {v_{\gamma}}\right )}{v_\gamma} +
\derII {\hat X} {\left(T_{\hat X \hat Y}(v_{\gamma})\right)} {v_\gamma} = $$
$$ \derII {\hat Y} {\left(\derII {\hat X} X {v_{\gamma}}\right )}{v_\gamma} +
R_{\hat X \hat Y} X (v_\gamma)  + \derII {\hat X} {\left(T_{\hat X \hat Y}(v_{\gamma})\right)} {v_\gamma}\tag 6.2$$
and this gives the claim since $\derII {\hat X} X {v_{\gamma^{(s)}}} = 
\derII {\doubledot \gamma} {\dot \gamma} {v_{\gamma^{(s)}}}
\equiv 0$ for any $s$.\par
The converse is proved  using suitable modifications of the
 arguments used for the analogous result 
in Riemannian or Hermitian geometry (see e.g. the proof of 
Prop.VII.1.1 in \cite{KN} vol. II) .\qed
\enddemo

\enddocument
\bye